\documentclass[12pt]{amsart}

\usepackage{amsmath}
\usepackage{amsthm}
\usepackage{amssymb}
\usepackage{amsfonts}
\usepackage{blkarray}
\usepackage{latexsym, graphicx}
\usepackage[top=1.5in, bottom=1.5in, left=1.5in, right=1.5in]{geometry}
\usepackage{color, enumitem}
\usepackage{tikz, xypic}
\usepackage{hyperref}
\usepackage[mathscr]{euscript}

\usetikzlibrary{arrows}
\usetikzlibrary{decorations.markings,decorations.pathmorphing,arrows}
\pgfarrowsdeclare{mytip}{mytip}{%
  \setlength{\arrowsize}{1pt}
  \addtolength{\arrowsize}{.5\pgflinewidth}
  \pgfarrowsrightextend{0}
  \pgfarrowsleftextend{-5\arrowsize}
}{%
  \setlength{\arrowsize}{1pt}
  \addtolength{\arrowsize}{.5\pgflinewidth}
  \pgfpathmoveto{\pgfpoint{-5\arrowsize}{1.5\arrowsize}}
  \pgfpathlineto{\pgfpointorigin}
  \pgfpathlineto{\pgfpoint{-5\arrowsize}{-1.5\arrowsize}}
  \pgfusepathqstroke
}
\tikzset{myarrow/.style={ decoration={bent,aspect=0.3, markings,mark=at
  position 0.5 with {\arrow[scale=1.2]{latex'}}}, postaction=decorate}}
\tikzset{myarrowshort/.style={ decoration={bent,aspect=0.3, markings,mark=at
  position 0.3 with {\arrow[scale=1.2]{latex'}}}, postaction=decorate}}
\tikzset{myarrowshorter/.style={ decoration={bent,aspect=0.3, markings,mark=at
  position 0.2 with {\arrow[scale=1.2]{latex'}}}, postaction=decorate}}
\tikzset{->-/.style={decoration={markings, mark=at position #1 with
  {\arrow{>}}},postaction={decorate}}}

\tikzset{my_dot/.style={fill, circle, inner sep=0pt,minimum size=1.5pt}}
\tikzset{my_node/.style={fill, circle, inner sep=0pt,minimum size=3pt}}
\tikzset{inv/.style={fill, circle, inner sep=0pt,minimum size=0pt}}

\newtheorem{thm}{Theorem}[section]

\newtheorem{fact}[thm]{Fact}

\newtheorem{Definition}[thm]{Definition}
\newenvironment{definition}
  {\vspace{.25cm} \begin{Definition}\rm}{\end{Definition}}

\newtheorem{Remark}[thm]{Remark}
\newenvironment{remark}
  {\vspace{.25cm} \begin{Remark}\rm}{\end{Remark}}

\newtheorem{Example}[thm]{Example}
\newenvironment{example}
  {\vspace{.25cm} \begin{Example}\rm}{\end{Example}}

\newtheorem{Exercise}[thm]{Exercise}
\newenvironment{exercise}
  {\vspace{.25cm} \begin{Exercise}\rm}{\end{Exercise}}

\newtheorem{Redefinition}[thm]{Redefinition}
\newenvironment{redefinition}
  {\vspace{.25cm} \begin{Redefinition}\rm}{\end{Redefinition}}

\newtheorem{Observation}[thm]{Observation}
\newenvironment{obs}
  {\vspace{.25cm} \begin{Observation}\rm}{\end{Observation}}

\newtheorem{Theorem}[thm]{Theorem}
\newenvironment{theorem}
  {\vspace{.25cm} \begin{Theorem}\rm}{\end{Theorem}}

\theoremstyle{remark}

\newcommand \un {\underline}
\newcommand \wt {\widetilde}
\newcommand \Sch {\mathrm{Sch}}
\newcommand \mmgn {\mathcal{M}_{g,n}}
\newcommand \omgn {\ov{\mathcal{M}}_{g,n}}
\newcommand \M {\mathcal{M}}
\newcommand \trop {\mathrm{trop}}
\newcommand \Mtrg {M_{g,n}^{\mathrm{trop}}}
\newcommand \mgn {M_{g,n}^{\mathrm{trop}}}

\newcommand \G {\mathbf{G}}  
\newcommand \D {\Delta}

\newcommand{\Gr}{\operatorname{Gr}}
\newcommand{\mcX}{\mathcal{X}}
\newcommand{\mfX}{\mathfrak{X}}\newcommand{\mfx}{\mathfrak{X}}
\newcommand{\mcO}{\mathcal{O}}
\newcommand{\mcU}{\mathcal{U}}

\newcommand \RR {\mathbb{R}}
\newcommand \R {\mathbb{R}}
\newcommand \ZZ {\mathbb{Z}}
\newcommand \Z {\mathbb{Z}}

\newcommand \ra {\rightarrow}

\newcommand \tr \textrm
\newcommand \eps \epsilon
\newcommand \sig \sigma

\newcommand{\val}{\operatorname{val}}

\newcommand{\Aut}{\operatorname{Aut}}

\newcommand \Af {\mathbb{A}}
 
\newcommand \ov {\overline}
\newcommand \PP {\mathbb{P}}
\newcommand \CC {\mathbb{C}}
\newcommand \QQ {\mathbb{Q}}

\newcommand \del {\partial}
\newcommand \col {\colon}
\newcommand {\set}[1]{\{1,\ldots,#1\}}

\newcommand{\laurent}{\mathbb{C}(\!(t)\!)}
\newcommand{\puiseux}{\mathbb{C}\{\!\{t\}\!\}}
\newcommand{\power}{\mathbb{C}[\![t]\!]}
\newcommand{\Spec}{\operatorname{Spec}}
\newcommand{\Trop}{\operatorname{Trop}}
\newcommand{\an}{\mathrm{an}}
\newcommand{\norm}[1]{|\!|#1|\!|}

\title{Lectures on tropical curves and their moduli spaces}
\author{Melody Chan}\address{Department of Mathematics, Brown University, Providence, RI 02912}\email{mtchan@math.brown.edu}

\begin{document}

\maketitle

\begin{abstract}
These are notes for a series of five lectures on ``Moduli and degenerations of algebraic curves via tropical geometry'' delivered at the \href{http://moduli2016.eventos.cimat.mx/node/304}{CIMPA-CIMAT-ICTP School on Moduli of Curves}, February 29-March 4, 2016 in Guanajuato, Mexico.
\end{abstract}

\tableofcontents

\section*{Introduction}

These are notes for five lectures on moduli and degenerations of algebraic curves via tropical geometry.  What do I mean by degenerations of algebraic curves?  The basic idea is that one can get information about the behavior of a smooth curve by studying one-parameter families of smooth curves, which degenerate in the limit to a singular curve, instead.  The singular curve typically has many irreducible components, giving rise to a rich combinatorial structure.  This technique obviously relies on having a robust notion of family of curves, that is, a moduli space.  Thus moduli spaces immediately come to the fore.

Tropical geometry is a modern degeneration technique.  You can think of it, to begin with, as a very drastic degeneration in which the limiting object is entirely combinatorial.  We will flesh out this picture over the course of the lectures.  It is also a developing field: exactly {\em what} tropical geometry encompasses is a work in progress, developing rapidly.

The powerful idea of using degenerations to study algebraic curves is at least several decades old and has already been very successful.  But recent developments in tropical geometry make it timely to return to and expand upon these ideas.  
The focus of these lectures will be on one beautiful recent meeting point of algebraic and tropical geometry: the {\em tropical moduli space of curves}, its relationship with the Deligne-Mumford compactification by stable curves, and its implications for the topology of $\mathcal{M}_{g,n}$.  I am not assuming any prior background in tropical geometry.

\bigskip

\section{From the beginning: tropical plane curves}\label{sec:embedded}

Let's start from scratch.  Centuries ago, some of the very first objects considered in algebraic geometry were just {\em plane curves}: the zero locus in $\PP^2$ of a homogeneous equation in three variables.  Only later was the perspective of studying curves abstractly, free from a particular embedding in projective space, developed.  We'll do the same, starting with some plane tropical curves, which are a special case of {\em embedded tropicalizations.}

\subsection{Embedded tropicalizations}
The natural setting for tropical geometry is over {\em nonarchimedean fields.}  Let $K$ be a field, and write $K^* = K\setminus\!\{0\}$ as usual.
\begin{definition}
A {\em nonarchimedean valuation} on $K$ is a map $v\colon K^*\ra \RR$ satisfying:
\begin{enumerate}
	\item $v(ab) = v(a)+v(b)$, and
	\item $v(a+b) \ge \min (v(a), v(b))$
\end{enumerate}
for all $a,b\in K^*$.  By convention we may extend $v$ to $K$ by declaring $v(0) = \infty$. 

The ring $R$ of elements with nonnegative valuation is the {\em valuation ring} of $K$.  Note that $R$ is a local ring; write $k = R/\mathfrak{m}$ for its residue field.
\end{definition}

\noindent If you are arithmetically minded, you might immediately think of the $p$-adic field $K=\QQ_p$, or $\CC_p$, the completion of its algebraic closure. For another example, take your favorite field, perhaps $\CC$, and equip it with the all-zero valuation, also known as the {\em trivial valuation}.  This example sounds unimportant, but it is theoretically important, because it permits a unified theory of tropicalization.

Another good example to keep in mind is $K=\laurent$, the field of {\em Laurent series,} with valuation $v(\sum_{i\in\ZZ}a_it^i) = \min \{i:a_i\ne 0\}$.  In this example, $K$ is discretely valued, in that $v(K ^* )\cong \ZZ$.  Note that the algebraic closure of $\laurent$ is the field of {\em Puiseux series} $\puiseux = \bigcup_n \CC(\!(t^{1/n})\!)$, whose elements are power series with bounded-denominator fractional exponents.  This field is a favorite of tropical geometers.

\begin{exercise}\label{ex:mintwice}
Let $a,b\in K$.  If $v(a)\ne v(b),$ then in fact $$v(a+b) = \min(v(a), v(b)).$$
In other words, for all $a,b\in K$, the minimum of $v(a), v(b), v(a+b)$ occurs at least twice.
\end{exercise}

\begin{exercise}
Suppose $K$ is a nonarchimedean field.  If $K$ is algebraically closed then its residue field $k$ is also algebraically closed.
\end{exercise}

\begin{remark}
Already, you can get a glimpse of {\em why} tropical geometry is performed over nonarchimedean fields.  For example, a curve $X$ over $K=\laurent$ can be regarded as a family of complex curves over an infinitesimal punctured complex disc $\Spec K$ around $t=0$.  The eventual {\em tropicalization} of this curve will be a metric graph (decorated with a little bit of extra stuff).  It records data about how the special fiber $t=0$ must be filled in according to the properness of $\ov\M_g$.  

\end{remark}

A {\bf valued extension} of valued fields is an extension $L/K$ in which the valuation on $L$ extends the valuation on $K$.  

\begin{definition} (Embedded tropicalization). \label{d:trop}
Fix a nonarchimedean field $K$.  Let $X$ be a subvariety of the algebraic torus $(\mathbb{G}_m)^n$.  So $X$ is defined by an ideal of the Laurent polynomial ring $K[x_1^\pm,\ldots,x_n^\pm]$.

The {\em tropicalization} of $X$ is the subset of $\RR^n$ 
$$\{ (v(x_1),\ldots,v(x_n)): (x_1,\ldots,x_n)\in X(L) \text{ for $L/K$ a valued extension}\}. $$
In particular, if $K$ is algebraically closed and nontrivially valued, then $\Trop(X)$ is the closure, in the usual topology on $\RR^n$, of the set
$$\{ (v(x_1),\ldots,v(x_n)): (x_1,\ldots,x_n)\in X(K) \}$$
of coordinatewise valuations of $K$-points of $X$.

\end{definition}

\noindent I am skimping a little bit on notation, and denoting all valuations by $v$, even over different fields. 

\bigskip

Let's immediately practice Definition~\ref{d:trop} in the case of a line.

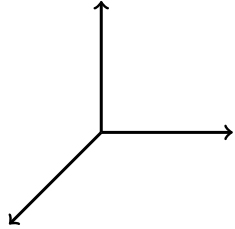
\begin{figure}
\begin{tikzpicture}[my_node/.style={fill, circle, inner sep=1.75pt}, scale=1.75]

\begin{scope}
\draw[line width = 0.04cm, ->] (0,0) to (1,0);
\draw[line width = 0.04cm, ->] (0,0) to (0,1);
\draw[line width = 0.04cm, ->] (0,0) to (-.7,-.7);
\end{scope}
\end{tikzpicture}

\caption{The tropical line from Example~\ref{ex:line}}
\label{fig:line}
\end{figure}

\begin{example}\label{ex:line}(A line in the plane).
Let $f(x,y)=x+y-1$ and let $X = V(f)\subset \mathbb{G}_m^2$.  So $X$ is $\PP^1$ minus 3 points.  What is $\Trop(X)$?

\bigskip

\noindent Answer: Suppose $x$ and $y$ are such that $x+y-1 = 0$. By Exercise~\ref{ex:mintwice},
$$\text{the minimum of $v(x), v(y),$ and $v(1)$ is attained at least twice.}$$
In other words, 
$$\Trop(X) \subseteq \{(z,w)\in \RR^2: \text{ the minimum of $z,w$, and $0$ occurs at least twice} \}.$$
You can draw this latter set: it is polyhedral, consisting of three rays from the origin in the directions of the standard basis vectors $e_1, e_2,$ as well as $-e_1-e_2$.  See Figure~\ref{fig:line}.  The content of the next theorem is that the containment above is actually an equality.
\end{example}

We state the theorem below in a baby case, the case of curves in the plane. But it holds verbatim for arbitrary hypersurfaces.  See \cite[Theorem 3.1.3]{ms} for the full statement and its proof.\footnote{The statement there also gives an equivalent formulation in terms of Gr\"obner initial ideals, which is key for computations, and which I won't talk about at all.}
\begin{thm}[Kapranov's Theorem]  Let $$f = \sum_{(i,j)\in \ZZ^2} \!\!c_{ij} x^iy^j \quad \in K[x^\pm,y^\pm].$$  Then 
$$\Trop(X) = \{(z,w)\in\RR^2: \min_{(i,j)\in\ZZ^2} v(c_{ij}) + iz + jw \text{ occurs at least twice.}\}$$
\end{thm}

\begin{remark}
Notice that the expression $\min \left( v(c_{ij}) + iz + jw\right)$ is obtained from $f$ by 
\begin{itemize}
\item replacing addition by minimum, 
\item replacing multiplication by $+$, and
\item replacing scalars $c_{ij}$ by their valuations.
\end{itemize}  This explains the slogan you may hear that {\em tropical geometry is the algebraic geometry of the min-plus semiring $(\RR\cup\{\infty\}, \min, +)$.}  
This may also help explain the naming of the field of tropical geometry, which was in honor of the Brazilian mathematician and computer scientist Imre Simon, a pioneer in the study of the min-plus semiring. \end{remark}

\begin{remark}
It is not nearly as straightforward to compute tropicalizations of subvarieties of $\mathbb{G}_m^n$ that are not hypersurfaces.  It can however be done by A.~N.~Jensen's software \texttt{gfan}, using Gr\"obner methods \cite{gfan}.  See the \texttt{gfan} manual for details.
\end{remark}

\begin{exercise}
What are all possible tropical lines in the plane, i.e., subsets of $\RR^2$ of the form $\Trop(X)$ where $X = V(ax + by + c)$?\end{exercise}

\medskip
\subsection{A very, very short treatment of Berkovich analytifications}  All of this could be said much more elegantly using the language of Berkovich spaces \cite{berk1}.  We now assume that $K$ is a {\em complete} valued field. This means that $K$ is complete as a metric space, with respect to its nonarchimedean valuation.\footnote{The valuation $v$ on $K$ defines a norm on $K$ by setting $|a| = \operatorname{exp}(-v(a)).$} Conceptually, the assumption that $K$ is complete is not such a big deal, because we can always base change from a given field $K$ to its completion $\hat{K}$.

\begin{example}
$\QQ_p$ is complete by construction.  Any trivially valued field is of course complete.  On the other hand, the field of Puiseux series $\puiseux$ is not complete. What is its completion? 
\end{example}

Let $X$ be finite type scheme over $K$.  We shall define the {\em Berkovich analytification} $X^{\an}$, a locally ringed topological space associated to $X$.  Actually, my plan is to entirely ignore the structure sheaf of analytic functions of $X^\an$.  So we will just regard $X^\an$ as a topological space for the duration.

I'll do things in an unconventional order, starting with a very quick way to say what the points of $X^\an$ are.  I find this simple description very useful---especially when $X$ is some kind of moduli space. 

\begin{definition}(Points of the Berkovich space).\label{d:berkpts}
Let $X$ be a finite type scheme over a complete nonarchimedean valued field $K$.  The points of the {\em Berkovich analytification} $X^\an$ are in bijection with maps $\Spec L \rightarrow X$ for all valued field extensions $L/K$, modulo identifying $\Spec L \ra X$ with $\Spec L' \ra \Spec L \ra X$, where $L'/L$ is again a valued field extension.
\end{definition}

Does that sound strange?   It should be compared with the more familiar situation of a scheme $Y$ over any field $K$, with no valuations in sight.  Then you can see for yourself the following way to name the points of $Y$: points of $Y$ correspond to maps $\Spec L \ra Y$ for all extensions $L/K$, modulo identifying $\Spec L \ra Y$ with $\Spec L' \ra \Spec L \ra Y$ for all further extensions $L'/L$.

Next we will define the topology on $X^\an$. 

\begin{definition}(Topology on the Berkovich space, affine case.)\label{d:berktop}
We continue to assume that $K$ is a complete nonarchimedean valued field.  Let $X = \Spec A$ be an affine scheme of finite type over $K$.  

We take $X^\an$ to have the coarsest topology such that for all $f\in A$, the function 
\begin{equation}\label{eq:nu}
\nu_f\colon X^\an \longrightarrow \RR
\end{equation} 
$$(\Spec L\stackrel{p}{\ra} X) \,\longmapsto\, v(p^\#f)$$ is continuous.
Here $v$ denotes the valuation on the valued field $L$ and $p^\#\colon A\rightarrow L$ is the map of rings coming from $p$.
\end{definition}

Now for an arbitrary finite type scheme $X$ over $K$ that is not necessarily affine, the topological space $X^\an$ is obtained by taking an affine open cover of $X$, analytifying everything separately, and then gluing.\footnote{Usually, the points of $(\Spec A)^\an$ are described as multiplicative seminorms $|\!|\cdot|\!|_p$ on $A$ extending the norm on $K$, equipped with the coarsest topology such that for every $f \in A$, the map $X^\an \rightarrow \RR$ sending $$\norm{\cdot}_p \mapsto \norm{f}_p$$ is continuous.
It's not hard to describe the correspondence between Definitions~\ref{d:berkpts} and this definition.  A seminorm $\norm{\cdot}_p$ corresponds to the map $$\Spec\, \operatorname{Frac} (A/\operatorname{ker}(\norm{\cdot}_p)) \ra \Spec A.$$}

Now let's try Definition~\ref{d:trop} over again:

\begin{redefinition} (Embedded tropicalization, again).\label{d:trop2}
Let $X \subseteq (\mathbb{G}_m)^n$, given explicitly as $X = \Spec K[x_1^\pm, \ldots, x_n^\pm]/I$.  The {\bf tropicalization} of $X$, denoted $\Trop(X)$, is the image of the map $X^\an\ra \RR^n$ that sends, for $p\colon \Spec L \ra X$ a point of $X^\an$,
$$p \quad \longmapsto \quad (\nu_{x_1}(p),\ldots,\nu_{x_n}(p)).$$
\end{redefinition}
The maps $\nu_{x_i}$ were defined in~\eqref{eq:nu}, and we set up the topology of $X^\an$ precisely so that each $\nu_{x_i}$ is continuous.  Thus
$\Trop(X)$ is, by Definition~\ref{d:berktop}, a {\em continuous} image of $X^\an$.  This is helpful!  For example, Berkovich tells us that $X^\an$ is connected if $X$ is connected \cite{berk1}.  Therefore, in this situation, $\Trop(X)$ is connected too.

We've just hinted at the fact that passing to analytifications can be a helpful perspective for viewing tropicalizations.  But actually, one could just as well say the reverse.  Namely, one of the reasons tropicalizations are useful is that they can provide a faithful ``snapshot'' of a piece of the much hairier\footnote{Almost literally.} and more complicated space $X^\an$.  See \cite{bpr,pay}.

\bigskip

\section{Abstract algebraic and tropical curves}

Next, what is an {\em abstract} tropical curve, and how does such a gadget arise from an algebraic curve over a valued field $K$?  That is the subject of this lecture.  The relationship with the previous lecture is as follows.  In the last section, we concerned ourselves with {\em embedded tropicalizations}, i.e.~tropicalizations of a subvariety of a torus or toric variety.  In this section, we are fast-forwarding many decades in the parallel story in the history of algebraic geometry, and treating curves now {\em in the abstract}, free from a particular embedding in projective space, say.  Also, just as in algebraic geometry, once this bifurcation between abstract and embedded tropicalization happens, it then becomes interesting to study the relationship between the two.  This is also a very interesting story (see  \cite{bpr, grw}), but I won't have time for it.  I also highly recommend \cite{bakerjensen} for a survey of the state of the art in tropical linear series and Brill-Noether theory.

Before launching into all the definitions, let me give one example in full. It's such a tiny example that I can guarantee that it's not that interesting on its own.  But it will serve as a little laboratory in which we can see all the definitions in action at once.

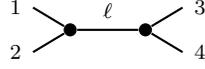
\begin{figure}
\vspace{.25cm}

\begin{tikzpicture}[my_node/.style={fill, circle, inner sep=1.75pt}, scale=2]

\begin{scope}
\node[my_node] (A) at (-.25,0){};
\node[my_node] (B) at (.25,0){};
\node[inv, label=right:$\scriptstyle 3$] (B1) at (.5, .15){};
\node[inv, label=right:$\scriptstyle 4$] (B2) at (.5, -.15){};
\node[inv, label=left:$\scriptstyle 1$] (A1) at (-.5, .15){};
\node[inv, label=left:$\scriptstyle 2$] (A2) at (-.5, -.15){};
\node[inv, label=above:$\scriptstyle\ell$] (L) at (0, 0){};
\draw[thick] (B) to (A);
\draw[thick] (B) to (B2);
\draw[thick] (B) to (B1);
\draw[thick] (A) to (A1);
\draw[thick] (A) to (A2);
\end{scope}
\end{tikzpicture}
\vspace{.25cm}

\caption{A picture of the abstract tropicalization of the curve $X/K$ of Example~\ref{ex:preview}.  We have drawn marked points as marked rays attached at vertices.}
\label{fig:M04}
\end{figure}

\begin{example}\label{ex:preview} (A preview). 
Let $K=\laurent$, with $R=\power$ its valuation ring. Let $\ell$ be a positive integer, and let $X/K$ be the projective plane curve with equation
\begin{equation}\label{eq:conic}
xy = t^\ell z^2.
\end{equation}
So, $X$ is just a smooth conic over $K$, but we regard $X$ as defining a germ of a family, with base parameter $t\ne 0$, of smooth plane conics in the complex projective plane.  Let's consider the four marked points 
$$p_1,p_2 = (\pm t^\ell: \pm 1: 1)\qquad p_3,p_4 = (\pm 1 : \pm t^\ell : 1)$$
on $X$.

Now, equation~\eqref{eq:conic} also defines a scheme $\mathfrak{X}/R$, in which the special fiber $\mathfrak{X}_k = \mathfrak{X} \times_R k$ has equation $xy=0$ in $\PP^2_\CC$.  That is, the special fiber is a union of two rational curves meeting at a node.  Furthermore, the horizontal closures $\ov{p_i}$ of the four marked points of $X$ do indeed meet $\mfX_k$ in four regular points, namely $(0:\pm1:1)$ and $(\pm1:0:1)$.  Note that $\mathfrak{X}$ will qualify as a {\em stable model} for $X$, as defined in Definition~\ref{def:smodel}.

The {\em abstract tropical curve} associated to $\mathfrak{X}$ will be the vertex-unweighted metric graph with two vertices and an edge between them of length $\ell$, with marked points $1,2$ on one vertex and $3,4$ on the other.  See Figure~\ref{fig:M04}.

\end{example}

\bigskip

%%%%%%%%%%%%%%%%%%%%%%%%%%%%%%%%%%%%

\bigskip

Now we'll launch into the definitions of stable curves, their dual graphs, and abstract tropicalizations.

\bigskip

\subsection{Stable curves}
Fix $k$ an algebraically closed field.  By a {\em curve} we shall mean a reduced, proper, connected scheme $X$ of dimension 1 over $k$.  The arithmetic genus of the curve is $h^1(X, \mcO_X).$  A {\em node} of $X$ is a point $p\in X(k)$ with the property that $\widehat{\mcO}_{X,p} \cong k[\![x,y]\!]/(xy)$.  A {\em nodal} curve is a curve whose only singularities, if any, are nodes.  

\begin{definition}(Stable $n$-pointed curves).
A nodal, $n$-marked curve of genus $g$ is $(X,p_1,\ldots,p_n)$, where $p_i \in X(k)$ are distinct nonsingular points of a genus $g$ nodal curve $X$.

We say that a nodal, marked curve $(X,p_1,\ldots,p_n)$ is {\bf stable} if $\Aut(X,p_1,$ $\ldots,p_n)$ is finite, that is, there are only finitely many automorphisms of the curve $X$ that fix each $p_1,\ldots,p_n$ pointwise.  This is often equivalently stated as follows: $(X,p_1,\ldots,p_n)$ is stable if the restriction of $\omega_X(p_1+\cdots+p_n)$ to every irreducible component of $X$ is a line bundle of positive degree.  Here $\omega_X$ denotes the dualizing sheaf of $X$.
\end{definition}

Notice that all smooth curves of genus $g\ge 2$ already have only finitely many automorphisms.  A smooth curve of genus $1$ has finitely many automorphisms once one fixes one marked point; and a smooth curve of genus $0$, also known as $\PP^1$, has finitely many automorphisms once one fixes three marked points.  So we could equally phrase the stability condition as follows:

\begin{obs}\label{obs:stability}
For every irreducible component $C$ of $X$, let $\phi\col C^\nu\ra C$ denote the normalization of $C$.  An $n$-marked nodal curve $(X,p_1,\ldots,p_n)$ is stable if and only if 
\begin{enumerate}
\item for every component $C$ of geometric genus 0, $$|C\cap \{p_1,\ldots,p_n\}| + |\{q\in C^\nu: \phi(q)\in X^ \mathrm{sing}\}| \ge 3;$$ 
\item for every component $C$ of geometric genus 1, $$|C\cap \{p_1,\ldots,p_n\}| + |\{q\in C^\nu: \phi(q)\in X^ \mathrm{sing}\}|  \ge 1.$$
\end{enumerate}
\end{obs}

\noindent (The second condition sounds misleadingly general. You can trace through the definition yourself to see that it excludes only one additional case, the case that the whole of $X$ is just a smooth curve of genus 1 with no marked points.)

\begin{exercise} \label{ex:2g-2+n}
Let $g, n\ge 0$.  Check that stable $n$-marked curves of genus $g$ exist if and only if $2g-2+n>0$.
\end{exercise}

\bigskip
\subsection{Stable models}

Let $K$ be an algebraically closed field that is complete with respect to a nonarchimedean valuation.  Good examples include the completion of the field of Puiseux series $\puiseux$ or the completion $\CC_p$ of the algebraic closure of the field $\QQ_p$.\footnote{By the way, you might complain that some of this theory can be developed with weaker hypotheses on the field $K$.  That is true.  For example, the stable reduction theorem holds for arbitrary complete nonarchimedean fields, up to {\em passing to a finite, separable field extension.}  I am taking this approach partly for expository ease, especially for one's first exposure to this material. It's kind of like learning algebraic geometry over $\CC$ first.  My other defense is that in the tropical context it is often not a big deal to pass to a possibly huge field extension, at least in theory.  See Definition~\ref{d:berkpts}, for example.}

As before, let $R$ denote the valuation ring of $K$ and let $k=R/\mathfrak{m}$ be the  residue field.  Recall that $\Spec R$ has two points $\eta$ and $s$, corresponding to the ideals $(0)$ and $\mathfrak{m}$ respectively.  If $\mfX$ is a scheme over $\Spec R$, then the {\em generic fiber} of $\mfX$ is the fiber over $\eta$; the {\em special fiber} is the fiber over $s$.

\begin{definition}(Models) \label{def:model}
If $X$ is any finite type scheme over $K$, then by a {\em model} for $X$ we mean a flat and finite type scheme $\mfx$ over $R$ whose generic fiber is isomorphic to $X$.
\end{definition}

Now let's define stable models.  First, let me forget about marked points, and just suppose that $X$ is a smooth curve over $K$.  

\begin{definition} (Stable models).\label{def:smodel}
Suppose $X$ is a smooth, proper, geometrically connected curve over $K$ of genus $g\ge 2$. 
A {\em stable model} for $X$ is a proper model $\mfx/R$ whose special fiber $\mfx_k =\mfx \times_R k$ is a stable curve over $k$.  
\end{definition}

\begin{definition} (Stable models, allowing marked points).\label{def:smodelmarked}
Now say $2g-2+n>0$ and suppose $(X,p_1,\ldots,p_n)$ is a smooth, $n$-marked, genus $g$ curve.  Then a stable model for $X$ is a proper model $\mfx/R$ with $n$ sections $\ov{p}_1,\ldots,\ov{p}_n\col \Spec R\ra \mfx$ restricting to the marked points $p_i$ on the general fiber, making the special fiber a stable $n$-marked curve of genus $g$ over $k.$  
\end{definition}

Let $2g-2+n>0$.  When does an $n$-marked genus $g$ curve $X/K$ admit a stable model? The answer is: {\em always}.  This is the content of the Stable Reduction Theorem of Deligne-Mumford-Knudsen, which also gives that the stable model is essentially unique.  More precisely, the version we are using here, for fields whose valuations are not necessarily discrete, goes back to \cite{bosch-lutkebohmert}; see also \cite{bpr, tyomkin, viviani}.\footnote{Again, the typical formulation of the stable reduction theorem says that if $X/K$ is a smooth curve, then there exists a finite separable field extension $K'/K$ such that $X\times_K K'$ admits a stable model.  Here, we've folded the need to pass to a finite field extension into the assumption that $K$ itself is algebraically closed.}    You can also see Harris and Morrison's book \cite[\S3.C]{harris-morrison} for a relatively explicit, algorithmic explanation of the Stable Reduction Theorem, at least in characteristic~0.

\bigskip

\subsection{Dual graphs of stable curves} \label{sec:dual graphs}
We are working towards the goal of associating a graph, with some vertex decorations and some edge lengths, to a smooth curve $X/K$.  The graph we are going to associate to $X$ is the {\em dual graph} of the special fiber of a stable model for $X$.  Basically, the dual graph of a stable curve $Y$ is a combinatorial gadget that records:
\begin{itemize}
\item how many irreducible components $Y$ has, and what their geometric genera are;
\item how the irreducible components of $Y$ intersect; and
\item the way in which the $n$ marked points are distributed on $Y$.
\end{itemize}

Now we will explain this completely, starting with the graph theory.  

\bigskip {\bf Conventions on graphs.}  All graphs will be finite and connected, with loops and parallel edges allowed.  (Graph theorists would call such objects finite, connected {\em multigraphs}.)  Remember that a graph $G$ consists of a set of vertices $V(G)$ and a set of edges $E(G)$.  Each edge is regarded as having two endpoints which are each identified with vertices of $G$, possibly the same. 

\begin{definition} \label{def:gmw} (Vertex-weighted marked graphs). 
A {\em vertex-weighted, $n$-marked  graph} is a triple $(G,m,w)$ where:
\begin{itemize}
\item $G$ is a graph;
\item $w\col V(G)\ra \ZZ_{\ge 0}$ is any function, called a {\em weight function}, and 
\item $m\col \{1,\ldots,n\}\ra V(G)$ is any function, called an $n$-marking.\footnote{Another common setup for marking a tropical curve is to attach infinite rays to a graph, labeled $\{1,\ldots,n\}$.  Our marking function $m$ is obviously combinatorially equivalent.}

\end{itemize}
The {\em genus} of $(G,m,w)$ is 
$$g(G) + \!\!\sum_{v\in V(G)} \!\!w(v)$$
where $$g(G) = |E|-|V|+1 $$ is the first Betti number of $G$, considered as a 1-dimensional CW complex, say.
\end{definition}

\begin{definition}\label{def:gmw-stable}(Stability for vertex-weighted marked graphs).
With $(G,m,w)$ as above, we'll say that $(G,m,w)$ is {\em stable} if for every $v\in V(G)$, 
$$2w(v) - 2 + \operatorname{val}(v) + |m^{-1}(v)| > 0.$$
Here $\operatorname{val}(v)$ denotes the graph-theoretic {\em valence} of the vertex $v$, which is defined as the number of half-edges incident to it.

\end{definition}

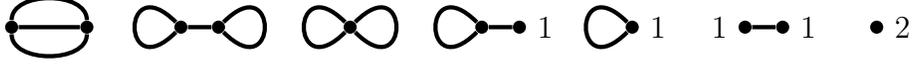
\begin{figure}\vspace{.5cm}

\begin{tikzpicture}[my_node/.style={fill, circle, inner sep=1.75pt}, scale=1]
% I
\begin{scope}
\node[my_node] (A) at (-.5,0){};
\node[my_node] (B) at (.5,0){};
\draw[ultra thick] (A)--(B);
\draw[ultra thick] (A) to [bend right=90] (B);
\draw[ultra thick] (A) to [bend left=90] (B);
\end{scope}
% II
\begin{scope}[shift = {(2,0)}]
\node[my_node] (A) at (-.25,0){};
\node[my_node] (B) at (.25,0){};
\draw[ultra thick] (A)--(B);
\draw[ultra thick] (A) to [out = 225, in = 135, looseness=20] (A);
\draw[ultra thick] (B) to [out = 45, in = 315, looseness=20] (B);
\end{scope}
% III
\begin{scope}[shift = {(4,0)}]
\node[my_node] (A) at (0,0){};
\draw[ultra thick] (A) to [out = 225, in = 135, looseness=20] (A);
\draw[ultra thick] (A) to [out = 45, in = 315, looseness=20] (A);
\end{scope}
% IV
\begin{scope}[shift = {(6,0)}]
\node[my_node] (A) at (-.25,0){};
\node[my_node, label=right:$1$] (B) at (.25,0){};
\draw[ultra thick] (A)--(B);
\draw[ultra thick] (A) to [out = 225, in = 135, looseness=20] (A);
\end{scope}
% V
\begin{scope}[shift = {(7.5,0)}]
\node[my_node, label=right:$1$] (A) at (.25,0){};
\draw[ultra thick] (A) to [out = 225, in = 135, looseness=20] (A);
\end{scope}
% VI
\begin{scope}[shift = {(9.5,0)}]
\node[my_node, label=left:$1$] (A) at (-.25,0){};
\node[my_node, label=right:$1$] (B) at (.25,0){};
\draw[ultra thick] (A)--(B);
\end{scope}
% VII
\begin{scope}[shift = {(11,0)}]
\node[my_node, label=right:$2$] (B) at (0,0){};
\end{scope}
\end{tikzpicture}
\vspace{.5cm}

\caption{The seven genus 2 stable vertex-weighted graphs with no marked points. The vertices have weight zero unless otherwise indicated.}
\bigskip

\bigskip
\label{f:genus2}
\end{figure}

Figure~\ref{f:genus2} shows the seven distinct stable vertex-weighted graphs of type $(g,n) = (2,0)$.  

\begin{exercise}\label{ex:m12}
Find the five stable, 2-marked weighted graphs of genus 1.  
\end{exercise}

\begin{exercise}\label{ex:finite-strata}
Prove that there are only finitely many stable marked, weighted graphs for fixed $g$ and $n$, up to isomorphism.  It may be helpful to consider the partial order of {\em contraction} that we will define in Section~\ref{sec:mtropg}.
\end{exercise}

\begin{definition}(Dual graph of a stable curve).  Let $k$ be an algebraically closed field, and let $(Y,p_1,\ldots,p_n)$ be a stable, $n$-marked curve over $k$.  

The {\em dual graph} of $(Y,p_1,\ldots,p_n)$ is the vertex-weighted, marked graph $(G,m,w)$ obtained as follows.

\begin{itemize}
\item The vertices $v_i$ of $G$ are in correspondence with the irreducible components $C_i$ of $Y$, with weights $w(v_i)$ recording the geometric genera of the components.
\item For every node $p$ of $Y$, say lying on components $C_i$ and $C_j$, there is an edge $e_p$ between $v_i$ and $v_j$.
\item The marking function $m\col \{1,\ldots,n\}\ra V(G)$ sends $j$ to the vertex of $G$ corresponding to the component of $Y$ supporting $p_j$.
\end{itemize}
\end{definition}

\noindent 
Note that by Observation~\ref{obs:stability}, $(G,m,w)$ is stable since $(Y,p_i)$ is stable.

\bigskip

\subsection{Abstract tropical curves}

A {\em metric graph} is a pair $(G,\ell)$, where $G$ is a graph, and $\ell$ is a function
$\ell\colon E(G) \rightarrow \R_{>0}$
on the edges of $G$. We imagine $\ell$ as recording real lengths on the edges of $G$.

An abstract tropical curve is just a vertex-weighted, marked {\em metric} graph:

\begin{definition}\cite{bmv, capsurvey, mik06} (Abstract tropical curves)
An {\bf abstract tropical curve} with $n$ marked points is a quadruple $\Gamma = (G,\ell,m,w)$ where:

\begin{itemize}
\item $G$ is a graph, 
\item $\ell\col E(G)\ra \RR_{>0}$ is any function, called a {\em length function}, on the edges,
\item $m\col \set{n} \ra V(G)$ is any function, called a {\em marking function}, and
\item $w\col V(G)\ra \ZZ_{\ge 0}$ is any function. 
\end{itemize}
The {\bf combinatorial type} of $\Gamma$ is the triple $(G,m,w)$, in other words, all of the data of $\Gamma$ except for the edge lengths.  We say that $\Gamma$ is {\bf stable} if its combinatorial type is stable.
The {\bf volume} of $\Gamma$ is the sum of its edge lengths.
\end{definition}
From now on, I will mean ``stable abstract tropical curve'' when I say ``abstract tropical curve,'' even if I forget to say so.

\bigskip

\noindent {\bf Hints of a tropical moduli space.}  Informally, we view a weight of $w(v)$ at a vertex $v$ as $w(v)$ loops, based at $v$, of infinitesimally small length. Each infinitesimal loop contributes 1 to the genus of $C$. 
Permitting vertex weights will ensure that the moduli space of tropical curves, once it is constructed, is complete. That is, a sequence of genus $g$ tropical curves obtained by sending the length of a loop to zero will still converge to a genus $g$ curve.

Of course, the real reason to permit vertex weights is so that the combinatorial types of genus $g$ tropical curves correspond precisely to dual graphs of stable curves in $\omgn$, and that the eventual moduli space will indeed be the boundary complex of $\M_{g,n}\subset \omgn$.

\medskip
\subsection{From algebraic to tropical curves: abstract tropicalization}

Now let's put everything together. 
We continue to let $K$ be an algebraically closed field, complete with respect to a nonarchimedean valuation.  Let $2g-2+n>0$.  Suppose $(X,p_1,\ldots,p_n)$ is a smooth, proper, $n$-marked curve over $K$ of genus $g$.  Let us extract a tropical curve from the data of $(X,p_i).$  

The procedure will go like this.  First we will extend $X$ to a family $\mfX$ over $\Spec R$ along with $n$ sections $\Spec R\ra \mfX$, filling in a stable, $n$-marked curve  of genus $g$ over $k$ in the special fiber.  The fact that this is possible is the {Stable Reduction Theorem}.

Then we will associate to $\mfX$ the {\em vertex-weighted dual graph} of the special fiber $\mfX_k$.  It only remains to equip the edges of the dual graph with real lengths.  We do this as follows: for every node $q$ of $\mfx_k$, say lying on components $C_i$ and $C_j$, the completion of the local ring $\mathcal{O}_{\mfx,q}$ is isomorphic to $R[\![x,y]\!]/(xy-\alpha)$ for some $\alpha\in R$, and $v(\alpha)>0$ is independent of all choices.  (See e.g~\cite{viviani}, and see Remark~\ref{rem:berkovich} below, for more on this independence.)  Then we put an edge $e_q$ between $v_i$ and $v_j$ of length $v(\alpha)$.  The result is a stable vertex-weighed, marked metric graph. 
See again Example~\ref{ex:preview}.  Summarizing:

\begin{definition}\label{def:abstract-trop}(Abstract tropicalization)
Let $K$ be an algebraically closed field, complete with respect to a nonarchimedean valuation. Suppose $(X,p_1,$ $\ldots, p_n)$ is a smooth, proper, $n$-marked curve over $K$ of genus $g$.  The {\bf abstract tropicalization} of $(X,p_i)$ is the dual graph of the special fiber of a stable model $(\mfX,\ov{p}_i)$ for $(X,p_i)$, declaring an edge corresponding to a node $q$ to have length $v(\alpha)$ if the local equation of $q$ in $\mfX$ is $xy-\alpha$.   
\end{definition}

\begin{remark}
We can now take Definition~\ref{def:abstract-trop} and extend it quite painlessly, to tropicalize {\em stable} curves, not just smooth ones.  In this situation, the local equation of a node in the special fiber may be of the form $xy = 0$; in other words, the node may have simply persisted from the general fiber.  Thus the natural result of tropicalization is a stable {\em extended} tropical curve: just like a tropical curve, but with edge lengths taking values in $\RR_{>0}\cup\{\infty\}$.  
\end{remark}

\begin{remark}\label{rem:berkovich}
Another way to say this whole story is that abstract tropicalization sends $X/K$ to its {\em Berkovich skeleton}: the minimal skeleton of $X^\an$, with respect to the $n$ marked points $p_1,\ldots,p_n$, equipped with the skeleton metric. See \cite[\S5]{bpr}.  (Actually, the Berkovich skeleton has $n$ infinite rays attached to the vertices, for the $n$ marked points).  There is a lot to be said here, but the main point at the moment is that the interpretation of the abstract tropicalization of $X$ as its Berkovich skeleton shows that it's canonically associated with $X$: the construction is in fact independent of all choices.
\end{remark}
\bigskip

\section{Definition of the moduli space of tropical curves}\label{sec:mtropg}
It is time to construct the moduli space of tropical curves.  This construction is due to Brannetti-Melo-Viviani \cite{bmv} and subsequently Caporaso \cite{capsurvey}, building on work of Mikhalkin \cite{mik06} and with antecedents in related constructions of Gathmann-Markwig \cite{gathmann-markwig, markwig-thesis}.  Actually, many of the ideas can be traced back even further to the work of Culler-Vogtmann \cite{cv}.

 Fix $g$ and $n$ with $2g-2+n>0$.  
Suppose we fix a single combinatorial type $(G,m,w)$ of type $(g,n)$, and allow the edge lengths $l$ to vary over all positive real numbers.  Then we clearly obtain all tropical curves of that type. This motivates our construction of the moduli space of tropical curves below. We will first group together curves of the same combinatorial type, obtaining one cell for each combinatorial type. Then, we will glue our cells appropriately to obtain the moduli space.    

To make this construction, for the moment we will just follow our noses combinatorially.  But the whole point of the next lectures will be that the space we get out the other side is a good one algebro-geometrically: it can be identified with the {\em boundary complex} of the Deligne-Mumford compactification $\omgn \supset \M_{g,n}$.

Let's begin. First, fix a combinatorial type $(G,m,w)$ of type $(g,n)$. What is a parameter space for all tropical curves of this type? 
Our first guess might be a positive orthant $\R_{>0}^{|E(G)|}$, that is, a choice of positive length for each edge of $G$.
But we have overcounted by symmetries of the combinatorial type $(G,m,w)$. For example, in the ``figure 8" depicted leftmost in Figure \ref{f:genus2}, the edge lengths $(2,5)$ and $(5,2)$ give the same tropical curve. 

Furthermore, with foresight, we will allow zero length edges as well, with the understanding that such a curve will soon be identified with one obtained by contracting those edges. This suggests the following definition:

\begin{definition}
  Given a combinatorial type  $(G,m,w)$, let
  the {\bf automorphism group} $\Aut(G,m,w)$ be the set of all permutations $\varphi: E(G) \to E(G)$ that arise from automorphisms of $G$ that preserve $m$ and $w$. 
    The group $\Aut(G,m,w)$ acts on the set $E(G)$, and hence on the 
  orthant $\R_{\geq 0}^{E(G)}$, with the latter action given by permuting coordinates. We define $\overline{C(G,m,w)}$ to be the quotient space
  \[ \overline{C(G,m,w)} = {\R_{\geq 0}^{E(G)}}/{\Aut(G,m,w)}. \]
\end{definition}

Next, we define an equivalence relation on the points in the union 
$$\coprod \overline{C(G,m,w)},$$ 
as $(G,m,w)$ ranges over all combinatorial types of type $(g,n)$. Regard a point $x \in \overline{C(G,m,w)}$ as an assignment of lengths to the edges of $G$. Now, given two points $x \in \overline{C(G,m,w)}$ and $x' \in \overline{C(G',m',w')}$, identify $x$ and $x'$ if one of them is obtained from the other by contracting all edges of length zero.  By {\em contraction}, we mean the following.  Contracting a loop, say based at vertex $v$, means deleting that loop and adding 1 to $w(v)$. Contracting a nonloop edge, say with endpoints $v_1$ and $v_2$, means deleting that edge and identifying $v_1$ and $v_2$ to obtain a new vertex whose weight we set to $w(v_1) + w(v_2)$. 

Let $\sim$ denote the equivalence relation generated by the identification we have just defined. Now we glue the cells $\overline{C(G,m,w)}$ along $\sim$ to obtain our moduli space:

\begin{definition}
  The {\bf moduli space} $\Mtrg$ is the topological space
  \[ \Mtrg \,\,:= \,\,\coprod \overline{C(G,m,w)} /\!\sim, \]
  where the disjoint union ranges over all combinatorial types of genus $(g,n)$, and $\sim$ is the equivalence relation defined above.
\end{definition}

\begin{figure}
\vspace{1cm}

\begin{tikzpicture}[my_node/.style={fill, circle, inner sep=1.75pt}, scale=1.25]
\begin{scope}[shift = {(-1.6,.8)}]
\node[my_node] (A) at (-.3,0){};
\node[my_node] (B) at (.3,0){};
\node[inv, label=left:{$\scriptstyle 1$}] (A2) at (-.6,-.0){};
\node[inv, label=right:$\scriptstyle 2$] (B2) at (.6,-.0){};
\draw[thick] (A) to [bend right=50] (B);
\draw[thick] (A) to [bend left=50] (B);
\draw[thick] (A) to (A2);
\draw[thick] (B) to (B2);
\end{scope}
\begin{scope}[shift = {(1.3,.8)}]
\node[my_node] (A) at (-.1,0){};
\node[my_node] (B) at (.25,0){};
\node[inv, label=right:$\scriptstyle 1$] (B1) at (.5, .25){};
\node[inv, label=right:$\scriptstyle 2$] (B2) at (.5, -.25){};
\draw[thick] (A) to [out = 210, in = 150, looseness=15] (A);
\draw[thick] (B) to (A);
\draw[thick] (B) to (B2);
\draw[thick] (B) to (B1);
\end{scope}
\begin{scope}[shift = {(0,-1)}]
\fill[gray!50] (0,0) -- (-1.6,0.8) -- (0,1.8) -- cycle;
\fill[gray!50] (0,0) -- (1.6,0.8) -- (0,1.8) -- cycle;
\node[inv] (V) at (0,0){};
\node[inv] (L2) at (-1.8,0.9){};
\node[inv] (M2) at (0,2){};
\node[inv] (R2) at (1.8,0.9){};
\draw[->, ultra thick] (V) --  (L2);
\draw[->, ultra thick] (V) --  (M2);
\draw[->, ultra thick] (V) --  (R2);
\draw[dotted, ultra thick] (V) to (-1,1.5);
\end{scope}
\end{tikzpicture}\vspace{1cm}

\caption{A picture of the tropical moduli space $M^{\trop}_{1,2}.$  Only the two top-dimensional strata are labeled.  We have drawn marked points as marked rays attached at vertices.}
\label{fig:M12}
\end{figure}
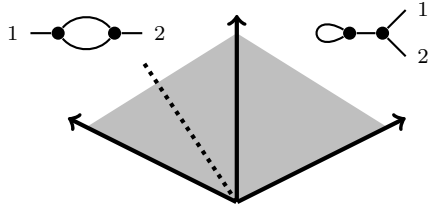

A picture of $M^{\trop}_{1,2}$ is shown in Figure~\ref{fig:M12}.  The picture is not entirely accurate, in that there is a 2-dimensional cone with a nontrivial symmetry which is drawn with a dotted line through it, which is supposed to remind us of the self-gluing of this cone induced by the symmetry.  

\begin{exercise}
Label the other cones of $M^{\trop}_{1,2}$ according to Exercise~\ref{ex:m12}.
\end{exercise}

\begin{exercise}
Verify that $\Mtrg$ is pure $(3g\!-\!3\!+\!n)$-dimensional, i.e., that the Euclidean dimension of every maximal cone $\overline{C(G,m,w)}$ of $\Mtrg$ is $3g\!-\!3\!+\!n$. 
\end{exercise}

Of course $\Mtrg$, being built out of cones, is contractible: e.g., it retracts onto its cone point, corresponding to the tropical curve, denoted $\bullet_{g,n}$, consisting of a single vertex with weight $g$, $n$ marked points, and no edges.  But the {\em link} of $\Mtrg$, meaning a cross-section of $\Mtrg$, is topologically very interesting and is a main character of these lectures.  

\begin{definition} \label{def:link} (Link of tropical moduli space)
The {\em link} $\D_{g,n}$ of $\Mtrg$ at the tropical curve $\bullet_{g,n}$ is the quotient of $\Mtrg \setminus \{\bullet_{g,n}\}$ induced by uniform scaling of edge lengths.  It can also be identified with the subspace of $\Mtrg$ parametrizing tropical curves of volume 1.\footnote{Having said all of that, you can find a slightly different way of describing $\Mtrg$, as a colimit of a diagram of rational polyhedral cones over the appropriate category of graphs, in \cite[\S4]{acp}, \cite[\S2]{cgp16}.  This colimit presentation equips $\Mtrg$ with the structure of a generalized cone complex, and its link has the structure of a smooth generalized $\Delta$-complex in the sense of \cite[\S3]{cgp16}.}

\end{definition}

%%%%%%%%%%%%%%%%%%%%%%%%%%%%%%%%%%%%

%%%%%%%%%%%%%%%%%%%%%%%%%%%%%%%%%%%%%%%%%%%%%%%%%%%%%%%%

\bigskip

\section{Boundary complexes of toroidal embeddings}

We are working towards the fundamental statement that {\em the link of the tropical moduli space of curves is the boundary complex of the stable curves compactification of $\M_{g,n}$.} This is one of the main results of \cite{acp}.  It is this identification that allows us to re-examine the boundary complex of $\mmgn\subset \ov{\mmgn}$ as a {worthy combinatorial moduli space} in its own right, and to obtain new results about the topology of $\mmgn$ using tropical geometry techniques. I will describe some of those applications in the last lecture.   Right now, we will answer the question: what is a {\em toroidal embedding}, and what is its {\em boundary complex}?

\subsection{\'Etale morphisms} Let's recall the definition of an \'etale morphism of schemes, which provides a more flexible notion of neighborhoods than do Zariski open neighborhoods.  Let $f\col X\ra Y$ be a morphism of schemes of finite type over a field $k$.  Then $f$ is called {\em unramified} if for all $x\in X$, letting $y= f(x)$, we have that $\mathfrak{m}_y\mathcal{O}_{X,x} = \mathfrak{m}_x$, and furthermore $k(y)$ is a separable field extension of $k(x)$.  Then $f$ is {\'etale} if it is both flat and unramified.  

\'Etale morphisms are, roughly, the algebro-geometric local analogue of finite covering spaces.  For example:

\begin{example}
Let $k=\CC$.  Then the map $\Af^1\ra \Af^1$ sending $z\mapsto z^n$ is \'etale away from $z=0$. 
\end{example}

\subsection{Toroidal embeddings}

The theory of toroidal embeddings is due to Kempf-Knudsen-Mumford-Saint-Donat \cite{kkmsd}.  Let $k$ be an algebraically closed field.  Let $U$ be an open subvariety of a normal variety $X$ over $k$.  We say that $U\subset X$ is a {\em toroidal embedding} if it is locally modeled by toric varieties.  More precisely, it is a toroidal embedding if for every $x\in X$, 
$$\widehat{\mathcal{O}}_{X,x} \cong \widehat{\mathcal{O}}_{Y_\sigma, y}$$
where $y\in Y_\sigma$ is a point in an affine toric variety $Y_\sigma$ with torus $T$, and furthermore, the ideals of $X\setminus U$ and $Y_\sigma\setminus T$ correspond in the respective completed local rings.

We'll write $D_1,\ldots,D_r$ for the irreducible components of $X-U$.  There are two cases, one of which makes the theory more intricate: we say that a toroidal embedding has {\em self-intersections} if the components $D_i$ are not all normal.  

It is actually not a problem if you aren't familiar with toric varieties, because the most relevant example for our purposes is one you definitely know: the usual embedding of the torus $\mathbb{G}_m^n$ into $\Af^n$. Note that the complement $\Af^n\setminus \mathbb{G}_m^n$ is the union of $n$ coordinate hyperplanes, intersecting transversely.  

Indeed, toroidal embeddings whose local toric charts are all affine spaces are called normal crossings divisors, and this is the case we'll be most interested~in.

\begin{definition}(Normal crossings and simple normal crossings)
Let $X$ be a normal variety, and $D$ a divisor. We say $D$ is a {\em normal crossings divisor} if  for every $x\in X$, we have $\widehat{\mathcal{O}}_{X,x} \cong k[\![x_1,\ldots,x_n]\!]$ and the equation of $D$ in $\widehat{\mathcal{O}}_{X,x}$ is $x_1\cdots x_i$ for some $i$.  Equivalently, $U=X-D \subset X$ is a toroidal embedding locally modeled by affine spaces.

We say $D$ is {\em simple normal crossings} if in addition $D$ has no self-intersections. 
\end{definition}

\begin{example}
The nodal cubic $V(y^2 = x^2+x^3)$ in $\Af^2$ is a normal crossings divisor, but not a simple normal crossings divisor.
\end{example}

\medskip
\subsection{Boundary complexes of toroidal embeddings}
The theory of boundary complexes for toroidal embeddings without self-intersection is due again to Kempf-Knudsen-Mumford-Saint-Donat \cite{kkmsd}.  For simplicity, we will state this theory in the case of simple normal crossings divisors, while emphasizing that both the work \cite{kkmsd} and the work of Thuillier \cite{thuillier} takes place in the more general case of toroidal embeddings.\footnote{Indeed, to a toroidal embedding $U\subset X$ without self-intersections, one associates a {\em rational polyhedral cone complex}, whose cones correspond to the toroidal strata of $U\subset X$ \cite{kkmsd}. Next, if $U\subset X$ is toroidal with self-intersections, one may associate a {\em generalized cone complex}, a more general object in which  self-gluings of cones are permitted  \cite{acp, thuillier}. 

How do these definitions specialize to Definitions~\ref{def:bdy} and~\ref{def:bdythuillier}, in the special cases of simple normal crossings and normal crossings, respectively?  In these cases, the cone complexes associated to $U\subset X$ are glued from {\em smooth} cones. The operation of replacing each smooth $d$-dimensional cone with a $(d\!-\!1)$-simplex produces the desired correspondence between these two pairs of definitions.  This is explained more precisely in~\cite[\S3]{cgp16}, in terms of an equivalence of categories between smooth generalized cone complexes and {\em generalized $\Delta$-complexes}.}

\begin{definition}\label{def:bdy}(Boundary complex, no self-intersections).
Suppose $U\subset X$ is an open inclusion whose boundary is simple normal crossings.  Let $D_1,\ldots,D_r$ be the irreducible components of $\del X = X-U$.  The {\em boundary complex} $\D(U\subset X)$, or just $\D(X)$, is the $\D$-complex on vertices $D_1,\ldots,D_r$ with a $d$-face for every irreducible component of an intersection $D_{i_1}\cap\cdots\cap D_{i_{d+1}}$. 
\end{definition}

\begin{example}
The boundary complex $\Delta(\mathbb{G}_m^n \subset \Af^n)$ is the simplex $\Delta^{n-1}$.  
\end{example}

Next, Thuillier recently extended the theory of boundary complexes in a way that is important for our applications, dropping the assumption that the $D_i$ are normal \cite{thuillier}.
\begin{definition}\label{def:bdythuillier}(Boundary complex, self-intersections).
Now let  $U\subset X$  have normal crossings boundary.  Let $V\ra X$ be an \'etale surjective morphism to $X$, such that $U_V = U\times_X V  \subset V$ is simple normal crossings, and let $V_2 = V\times_X V$, with $U_2 = U_V\times_X U_V$.  The boundary complex of $U\subset X$ is the coequalizer, in the category of topological spaces, of the diagram
$$\D(U_2 \subset V_2) \rightrightarrows \D(U_V \subset V).\footnote{Thuillier actually shows that this construction is independent of all choices, because in fact it is intrinsic to the Berkovich analytification of the pair $U\subset X$.  See \cite{thuillier} for the precise description.}$$
\end{definition}

\begin{figure}
\includegraphics[width=1.5in]{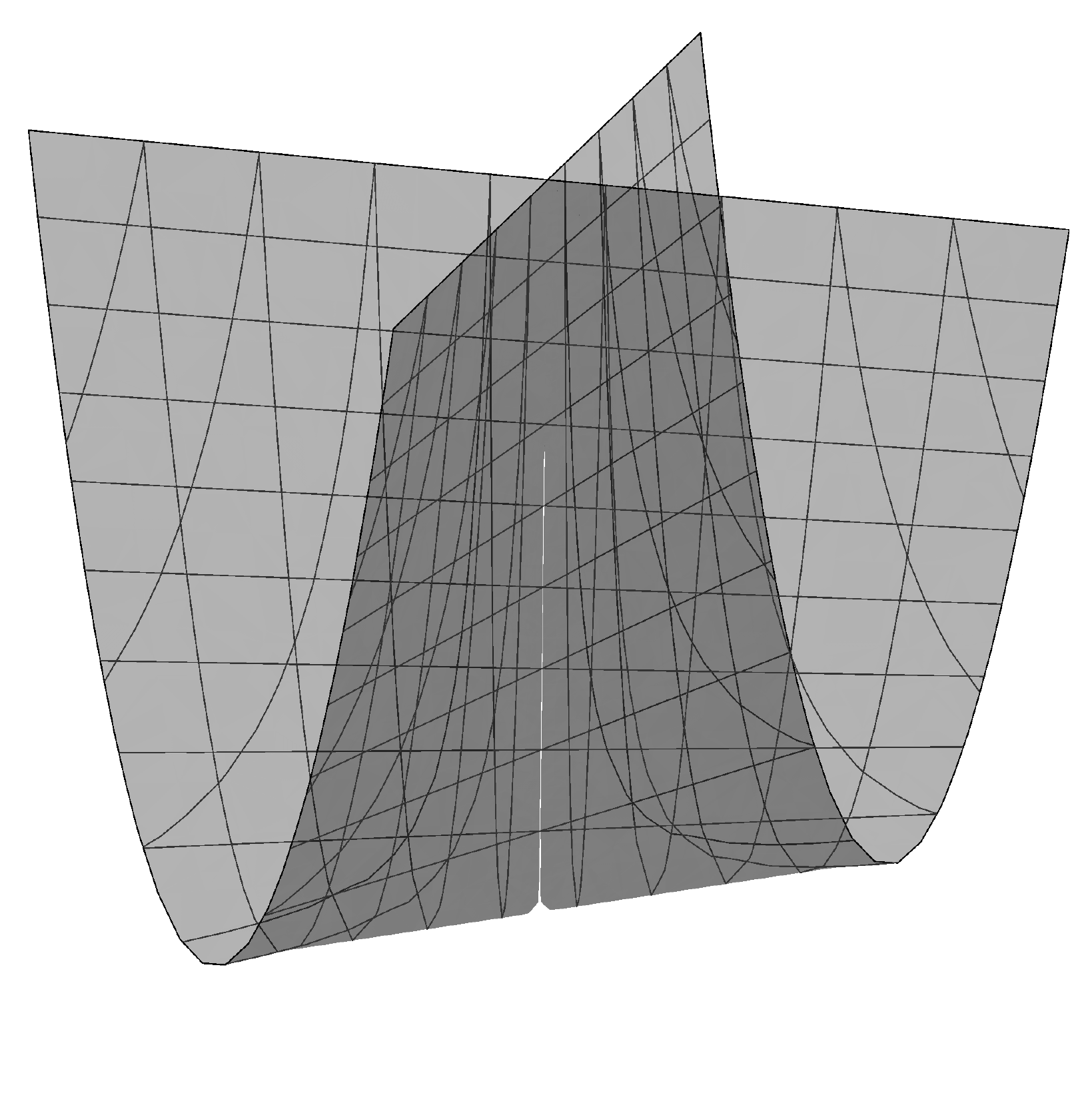}
\vspace{-.5cm}

\caption{The Whitney umbrella of Example~\ref{ex:whitney}.}
\label{fig:whitney}
\end{figure}

\begin{example}\label{ex:whitney} Let $k=\CC$.  Consider, as in \cite[Example 6.1.7]{acp}, the {\em Whitney umbrella} $$D = \{x^2y=z^2\} \quad \subset \quad X=\mathbb{A}^3\setminus \{y=0\},$$ 
drawn in Figure~\ref{fig:whitney}.
Let $U = X-D$.  We will explain why $\Delta(U\subset X)$ is a ``half-segment,'' meaning the quotient of a line segment  by a $\ZZ/2\ZZ$ reflection.\footnote{Of course, a line segment modulo a reflection is, topologically, just another segment.  There is a more abstract definition of a boundary complex in which the half-segment and segment are nonisomorphic, and only their {\em geometric realizations} as topological spaces are homeomorphic.  See \cite[\S3]{cgp16}.
}

Let $V\cong \mathbb{A}^2\times \mathbb{G}_m \rightarrow X$ be the \'etale cover of degree 2 given by a base change $y = u^2$.  Then $$D_V = D\times_X V = \{x^2u^2 - z^2 = 0\}$$ is simple normal crossings, and $D_2 = D_V\times_X D_V \cong D_V\times \Z/2\Z$, since $D_V$ is degree 2 over $D$.  Explicitly, one component of $D_2$ parametrizes pairs $(p,p)$ of points in $D_V$, and the other parametrizes pairs $(p,q)$ with $p\ne q$ lying over the same point of $D$. 

That means that $\Delta(V)$ is a segment and $\Delta(V_2)$ is two segments, and the two maps $\Delta(V_2)\rightrightarrows \Delta(V)$ differ by one flip.  So the coequalizer is a segment modulo a flip.

What is going on complex-analytically?  Let $Y =\{(0,y,0):y\ne 0\}$ be the umbrella pole.  It is a punctured affine line $\Af^1-\{0\}$ over $\CC$, to be visualized complex-analytically as a punctured plane.  The points $(0,u,0)$ and $(0,-u,0)$ in $D_0$ correspond to the two analytic branches of $D$ along $Y$ at the point $(0,y,0)$, where $y=u^2$.  The equations of the branches are $z=xu$ and $z=-xu$.  So taking $y$ around a loop around the puncture precisely interchanges the branches.  
\end{example}

\begin{exercise}
Compute the boundary complex of the complement of the nodal cubic $V(y^2 = x^2+x^3)$ in $\Af^2.$\end{exercise}

\bigskip

\section{Toroidal Deligne-Mumford stacks}

Boundary complexes can be defined for toroidal Deligne-Mumford stacks as well by following Thuillier's construction \cite{acp}.  The punchline will be that toroidal Deligne-Mumford stacks admit \'etale covers by toroidal schemes, and so Definition~\ref{def:bdythuillier} can be repeated with respect to toroidal \'etale covers with no changes.  
Here we'll give a quick-start guide to toroidal DM stacks.   We will have to skip many details in order to get anywhere in the allotted time. But I'll try to indicate exactly what I'm skipping.

\bigskip

\subsection{Categories fibered in groupoids}
Let us fix an algebraically closed field $k$, without the structure of a valuation.  
Let $\Sch_k$ denote the category of schemes over $k$.  Let me recall the following ``negative result'': it's impossible to define a scheme that deserves to be called a {\em fine moduli space} for genus $g$, $n$-marked curves. By a {\em fine moduli space} I mean a space $\mmgn$ such that maps $S$ to $\mmgn$ correspond, functorially, to families of genus $g$, $n$-marked curves over $S$.  The obstruction is that some curves have nontrivial automorphisms.  

Nevertheless, we can start by axiomatizing the desired property of a fine moduli space into a category whose objects are families of genus $g$, $n$-marked smooth, respectively stable, curves.  More precisely:

\begin{definition} \label{def:mmgn}(The category $\mmgn$). 
We denote by $\mmgn$ the category whose objects  are flat, proper morphisms $f\col X\ra B$ of $k$-schemes, together with $n$ sections $p_1,\ldots,p_n\col B\ra X$, such that the geometric fibers, with their $n$ marked points induced by the $p_i$, are smooth curves of type $(g,n)$.  The morphisms in $\mmgn$ are Cartesian diagrams
\begin{equation*}
      \xymatrix{
        X' \ar[r]\ar[d] & X \ar[d]\\
        B' \ar[r] & B.  }
  \end{equation*}
\end{definition}

\begin{definition} \label{def:omgn}(The category $\omgn$). 
The definition of $\omgn$ is the same as above, but with families of {\em stable} curves instead.
\end{definition}

These two categories come equipped with obvious functors to $\Sch_k$: take a family of curves $X\ra B$, and remember only its base $B$.  In fact, both $\mmgn$ and $\omgn$, along with their functors to $\Sch_k$, are examples of {\em categories fibered in groupoids}, or CFGs for short.  I won't state the condition that a category is {\em fibered in groupoids} over $\Sch_k$, but I will state that in our case it is the condition that pullbacks of families of curves exist and are unique up to unique isomorphism.

Here is another CFG, verifying that CFGs encompass $k$-schemes:

\begin{definition}\label{def:sch_}(The category $\un{S}$).  Let $S$ be any $k$-scheme.  The objects of the category $\un{S}$ are morphisms $X\ra S$ of $k$-schemes. The morphisms in $\un{S}$ are commuting triangles $X'\ra X\ra S$.  The functor $\un{S}\ra \Sch_k$ sends $(X\ra S)$ to $X$.
\end{definition}
\noindent In addition, the category $\un{S}$ determines the scheme $S$, in the sense made precise by Yoneda's Lemma.  In this case, the fact that $\un{S}$ is a CFG boils down to the fact that the  {\em composition} of two morphisms $X'\ra X\ra S$ exists and is unique---which is obvious.  

A {\em morphism} of CFGs $\mathcal{C}$ and $\mathcal{D}$ is just what you think: it is a functor $F\col \mathcal{C}\ra\mathcal{D}$ making a commuting triangle with the functors $\mathcal{C}\ra \Sch_k$ and $\mathcal{D}\ra \Sch_k$. 
Using the Yoneda correspondence, you can check (and make more precise): 
\begin{exercise}
To give a morphism $\un{S}\ra \mmgn$ is precisely to give a family of genus $g$, $n$-marked curves over $S$.
\end{exercise}

\medskip

\subsection{Fast forward: Deligne-Mumford stacks, and toroidal embeddings}
Now let's chat a little about Deligne-Mumford stacks.  Not all categories fibered in groupoids are schemes. (In other words, not all CFGs are of the form $\un{S}$ for some $k$-scheme $S$.) {Stacks}, and, even more restrictively, {\em Deligne-Mumford stacks}, are CFGs satisfying some extra conditions that make them behave a little more geometrically, even if they aren't exactly schemes. 

One of these requirements is the following, which we state as a fact about Deligne-Mumford stacks:

\begin{fact}\label{fact:etale}
If $\M$ is a Deligne-Mumford stack, then there is an \'etale, surjective morphism from a scheme $U$ to $\M$.\footnote{What do I even mean by saying that the morphism  $U\ra \M$ has a property like, say, \'etale?  The other condition for a stack $\M$ to be Deligne-Mumford is a representability condition: it says that for any morphisms $f\col S\ra \M$ and $g\col T\ra \M$ from a scheme, the {\em fiber product} $S\times_\M T$ (which I will not define) is again a scheme.  Then for any property $P$ of morphisms that is preserved by base change, we say that $f$ has property $P$ if $S\times_\M T \ra T$ is a map of schemes with property $P$.} 
\end{fact}

\noindent A very rough rephrasing is that locally, everywhere in a Deligne-Mumford stack there is a ``scheme covering space.''

Given Fact~\ref{fact:etale}, we can define toroidal Deligne-Mumford stacks by looking on \'etale atlases, and define their boundary complexes in exactly the same way as in Definition~\ref{def:bdythuillier}.
We continue to let $k$ be an algebraically closed field without valuation.  All our stacks are separated and connected over $k$.  

\begin{definition}\label{def:toroidaldm} (Toroidal Deligne-Mumford stacks).
An open substack $\mcU \subset \mcX$ of a Deligne-Mumford stack $\mcX$ is {\em toroidal} if for every \'etale morphism $V\ra \mcX$ from a scheme, the induced map of schemes $\mcU_V := \mcU\times_\mcX V \ra V$ is a toroidal embedding (of schemes).
\end{definition}

\begin{definition}\label{def:bdystack} (Boundary complexes of toroidal Deligne-Mumford stacks).
This is a reprise of Definition~\ref{def:bdythuillier}. Let $\mcU \subset \mcX$ be a toroidal Deligne-Mumford stack.  Let $V\ra \mcX$ be an \'etale cover by a scheme such that $\mcU_V \ra V$ is a toroidal embedding of schemes without self-intersections.  Then the {\em boundary complex} of $\mcX$ is the coequalizer, in the category of topological spaces, of
$$\D(V\times_\mcX V)\rightrightarrows \D(V).$$
\end{definition}

That's it. Abramovich-Caporaso-Payne show that Thuillier's work can be extended to the setting of DM stacks. In particular, when $\mcX$ is proper, the boundary complex of $\mcU\subset \mcX$  can be found intrinsically inside the Berkovich analytification of the coarse moduli space of $\mcX$.  \footnote{In fact, it can be found intrinsically inside the analytification of the stack $\mcX$ itself \cite[IV]{ulirsch-thesis}.}

\bigskip

\section{Compactification of $\M_{g,n}$ by stable curves}

We continue to let $k$ be an algebraically closed field, with no valuation.
At this point, I want to recall some of the essential facts about the Deligne-Mumford-Knudsen moduli stacks $\M_{g,n}$ and $\omgn$ \cite{deligne-mumford, knudsen1, knudsen2, knudsen3}.  

\begin{fact} The category $\omgn$ of smooth genus $g$, $n$-marked curves over $k$, defined in Definition~\ref{def:omgn}, is a smooth, proper Deligne-Mumford stack containing  $\mmgn$ an open substack. The inclusion $\mmgn\subset \omgn$ is toroidal, indeed normal crossings (though far from simple normal crossings).
\end{fact}

\medskip

\subsection{The boundary strata $\M_\G$ of $\omgn$}  
Moreover, the strata of the boundary $\omgn\setminus\mmgn$ are naturally indexed by genus $g$, $n$-marked combinatorial types $\G=(G,m,w)$, according to the dual graphs of the stable curves that they parametrize. I would now like to describe these strata, which we'll denote $\M_\G$.  This description follows \cite[\S3.4]{acp} and the correctness of this description is proved in \cite[\S12.10]{acg2}.

Fix a combinatorial type $\G=(G,m,w)$.  For each vertex $v$, let $n_v = \val(v) + |m^{-1}(v)|$ where $\val(v)$ is the valence of $v$.  Let
$$\wt{\M_\G} = \prod_{v\in V(G)} \M_{w(v), n_v}.$$
If you think about it, $\wt{\M_\G}$ can be identified with the moduli space of $n$-marked genus $g$ stable curves, {\em together with a chosen isomorphism of the dual graph with $\G$.}  To get rid of that choice of isomorphism, we take the stack quotient $[\wt{\M_\G}/\Aut(\G)]$.   The theorem is then that 
there is a canonical isomorphism $$\M_\G \cong [\wt{\M_\G}/\Aut(\G)].$$
An explicit example is given in Example~\ref{ex:M13} below.

(It is worth noting that this stratification is {\em inclusion-reversing} with respect to the corresponding stratification of $\mgn$ by combinatorial type. The more edges there are in $\G$, the smaller the stratum $\M_\G$ is, and the larger $\ov{C}(\G)$ is in the tropical moduli space.)

Now, given a point $p\in \M_\G$ corresponding to a stable curve $C$, we may describe an \'etale neighborhood $V_p$ of $p$ in $\omgn$ in which the boundary can be identified with the $d$ coordinate hyperplanes inside $\Af^d$.  The boundary of $\omgn$ in this neighborhood $V_p$ is a union of irreducible divisors $D$ with simple normal crossings; each $D_i$ corresponds to an edge of $\G$ and parametrizes local smoothings of the corresponding node.  So the boundary complex of $V_p$ is just a simplex $\D^{E(\G)-1}$.  

But there is in fact {\em monodromy} manifested in the coequalizer $$\D(V_p\times V_p)\rightrightarrows \D(V_p),$$ and it turns out that this monodromy identifies the coequalizer of the diagram above with $$\D^{E(\G)-1}/\!\Aut(\G),$$ where $\Aut(\G)$ acts by permutation on $E(\G)$.  So there is concordance on the level of strata with $\mgn$!

\begin{example}\label{ex:M13}
Let's see everything at work in the following specific example of a stratum in $\overline{\M}_{1,3}$.  

Let $\G$ be the combinatorial type below.

\begin{figure}[h!]
%%%%%%%%%%%%%%%%%%%%%%%%%%%%%%%%%%%%
\begin{tikzpicture}[my_node/.style={fill, circle, inner sep=1.75pt}, scale=1.25]
\begin{scope}[shift = {(-1.6,.8)}]
\node[my_node] (A) at (-.3,0){};
\node[my_node] (B) at (.3,0){};
\node[inv, label=left:{$\scriptstyle 1$}] (A2) at (-.6,.2){};
\node[inv, label=left:{$\scriptstyle 2$}] (A3) at (-.6,-.2){};
\node[inv, label=right:$\scriptstyle 3$] (B2) at (.6,-.0){};
\draw[thick] (A) to [bend right=50] (B);
\draw[thick] (A) to [bend left=50] (B);
\draw[thick] (A) to (A2);
\draw[thick] (A) to (A3);
\draw[thick] (B) to (B2);
\end{scope}
\end{tikzpicture}
\end{figure}
\noindent Consider the boundary stratum $\M_\G$ of $\overline{\M}_{1,3}$.  Locally, it is a self-intersection of the boundary component whose dual graph is obtained from $\G$ by contracting either edge.

Let's describe $\M_\G$.  I'll assume $\operatorname{char}  k \ne 2$ in this example.  According to the discussion above we have $\wt{\M_\G} \cong \M_{0,4}$. Essentially, to give a stable curve $C$ with dual graph $\G$ {\em along with a fixed identification of the two nodes of $C$ with the two edges of $\G$}, we choose (up to projective equivalence) four distinct points $p_1,p_2,q_1,q_2$ on a $\PP^1$, with the understanding that $p_1$ will be marked 1, $p_2$ marked 2, and $q_1$ and $q_2$ will be the two points of attachment of the other rational curve.  Of course $\M_{0,4}$ is an honest variety: for example, fixing $p_1 = 0, p_2 = 1,$ and $q_1 = \infty$ identifies $\M_{0,4}$ with $\Af^1-\{0,1\}$.

Now $\M_\G$ is then the stack quotient $[\M_{0,4}/(\ZZ/2\ZZ)]$, where the action is the one that exchanges $q_1$ and $q_2$.  You can work out that with the identification $\wt{\M_\G}=\Af^1-\{0,1\}$ above, the action sends $a$ to $1-a$.  

Thus the quotient $\M_\G$ is a once-punctured plane with a $\ZZ/2\ZZ$-stacky point, corresponding to the fixed point $(0,1,\infty,1/2)$ of $\M_{0,4}$ under $\ZZ/2\ZZ$.  It is the stacky point that produces monodromy: walking around it interchanges the analytic branches of the boundary divisor that meet along it.  This example is just like Example~\ref{ex:whitney}, except that the punctured complex plane in that example  is now filled in with a $\ZZ/2\ZZ$-stacky point, and there is another (inconsequential) puncture elsewhere.  

The result is that for a point $p$ in $\M_\G$, the neighborhood $V_p$ has boundary complex a segment modulo a flip, just as in Example~\ref{ex:whitney}. And this is indeed a slice of the cell $\ov{C(\G)}$ corresponding to $\G$ in the tropical moduli space $M^{\mathrm{trop}}_{1,3}$.
\end{example}

I omitted many details here, but this can all be patched together to show:

\begin{thm}\cite{acp}\label{thm:acp}
There is a canonical identification of the link $\D_{g,n}$ of $\mgn$ with the boundary complex of the toroidal embedding $\mmgn\subset \omgn$.  
\end{thm}

\bigskip

\subsection{Applications to the topology of $\mmgn$}  What good is all of that? Here is one application to the cohomology of $\mmgn$.

Suppose $U$ is a smooth variety over $k$, or even a smooth Deligne-Mumford stack.  Let $ X$ be a normal crossings compactification of $U$.  It is a fact that the homotopy type of the boundary complex $\D(X)$ is independent of the choice of compactification $X$; that is, any other normal crossings compactification produces a homotopy equivalent complex \cite{danilov}.  This means that all topological invariants of $\D(X)$ are actually invariants of $U$ itself.  
One very interesting such invariant is the {\em rational homology} of $\D(X)$.  The reason that it is particularly interesting is as follows.

Set $k=\CC$.  Then there is a {\em weight filtration}, due to Deligne, on the cohomology of $U$
$$W_0 H^k(U,\QQ) \subset \cdots \subset W_{2k}H^k(U,\QQ) = H^k(U,\QQ).$$
Let's write  $\Gr^W_k \!H^j$ for the quotient $W_k H^j / W_{k-1} H^j$. Letting $d=\dim U$, I'll refer to $\Gr^W_{2d} \!H^*$ as the {\em top-weight cohomology}, since cohomology never appears in weights above $2d$.  The point is that there is a canonical identification
\begin{equation}\label{eq:id}
\widetilde{H}_{i-1}(\Delta(U\subset X), \QQ) \cong \Gr^W_{2d}H^{2d-i}(U,\QQ)
\end{equation}
of reduced, rational homology of the boundary complex, up to shifting degrees, with the top-weight rational cohomology of $U$. 
These facts all follow from Deligne's work on mixed Hodge structures \cite{deligne3} in the case of varieties; in the case of stacks, they can be proved following Deligne's ideas as in \cite[Appendix]{cgp16}.  In short: the top-weight slice of cohomology is combinatorially encoded in the boundary complex of any normal crossings compactification.

The identification~\eqref{eq:id}, along with Theorem~\ref{thm:acp}, allows us to study the rational cohomology of $\mmgn$ appearing in top weight exactly by studying the reduced rational homology of $\mgn$.  This is a useful shift in perspective, because it allows arguments from tropical geometry and metric graph theory to be employed to study these complexes. Along these lines, S.~Galatius, S.~Payne, and I have shown some general results on $\Delta(\M_{g,n})$, including a lower bound on the connectivity of those spaces.\footnote{We say that a space is $n$-connected if the homotopy groups $\pi_1,\ldots,\pi_n$ all vanish.}
When $g=1$, our results allow us to describe the whole situation pretty thoroughly:

\begin{theorem}\cite{cgp16}\label{thm:cgp}
\begin{enumerate}
\item
$\Delta(\M_{1,n})$ is homotopy equivalent to a wedge of $(n-1)!/2$ top-dimensional spheres, for $n\ge 3$.  (It is contractible when $n=1$ and $2$.) Therefore:
\item  For each $n\ge 1$, the top weight cohomology of $\M_{1,n}$ is
 \[
\Gr_{2n}^W H^i(\M_{1,n}, \QQ) \cong \left \{ \begin{array}{ll} \QQ^{(n-1)!/2}  & \mbox{ \ \ for $n \geq 3$ and $i = n$,} \\  
                                                 0 & \mbox{ \ \ otherwise.} \end{array} \right.
 \]
 Moreover, for each $n\ge 3$, the representation of $S_n$ on $\Gr^W_{2n} H^n(\mathcal{M}_{1,n}, \QQ)$ induced by permuting marked points can be described explicitly, as in \cite{cgp16}.

\item When $n\ge 3$, a dual basis for $\Gr^W_{2n}(\M_{1,n},\QQ)$ is given by the {\em torus classes} associated to the $(n-1)!/2$ terminal curves whose dual graph is a loop of $n$ once-marked $\PP^1$s.
\end{enumerate}
\end{theorem}
\noindent I should remark that when $n\ge5$, Theorem~\ref{thm:cgp}(2) also follows in principle from an earlier calculation by  E.~Getzler \cite{getzler}, as we explain in \cite{cgp16}.
There is also a companion paper for the case $g=2$ \cite{cha14}, in which I again use tropical techniques to show vanishing of {\em integral} homology of $\D(\M_{2,n})$ outside the top two degrees, and compute the top-weight Euler characteristic of $\M_{2,n}$ for every $n$.  Furthermore, using tropical geometry and {\em using a computer} one can fully compute the top-weight $\QQ$-cohomology of $\mathcal{M}_{g,n}$ in a range of cases, as we discuss in \cite{cha14, cgp16}.   I reproduce the computations for $\M_{2,n}$ from \cite{cha14} for your curiosity: the cohomology is concentrated in degrees $n+3$ and $n+4$ with ranks given in the table below.
\begin{table}[h]
\begin{tabular}{c|ccccccccc}
    \hline
   $n$ &        $0$ & $1$& $2$ & $3$ & $4$ & $5$  &  $6$ & $7$ & $8$\\
  \hline
  $\dim \Gr^W_{2d} H^{n+3}(\M_{2,n},\QQ)$ &  $0$ & $0$ & $1$ & $0$ & $3$ & $15$ & $86$ & $575$ & $4426$\\
  \hline
  $\dim \Gr^W_{2d} H^{n+4}(\M_{2,n},\QQ)$ &  $0$ & $0$ &$0$ & $0$ & $1$ & $5$ &  $26$ & $155$ &$1066$ \\
  \hline
 \end{tabular}
 \medskip
  \label{table:upto8}
\end{table}

\bigskip

\noindent {\bf Acknowledgments.}  Thank you very much to the organizers of the Moduli of Curves School at CIMAT for inviting me to lecture and to write these notes, and to Dan Abramovich, Ethan Cotterill, Sam Payne, and the anonymous reviewer for giving me comments on them.  I'm grateful to Dan Abramovich, Matt Baker, Lucia Caporaso, Joe Harris, Diane Maclagan, Sam Payne, and Bernd Sturmfels, and many others for many helpful conversations over the years.


\begin{thebibliography}{99}

\bibitem{acp} D.~Abramovich, L.~Caporaso, S.~Payne, The tropicalization of the moduli space of curves,  Ann.~Sci.~Ec.~Norm.~Sup\'er., 48 (2015), no. 4, 765--809.

\bibitem{acg2} E.~Arbarello, M.~Cornalba, P.~Griffiths, Geometry of algebraic curves. Volume II. 
With a contribution by Joseph Daniel Harris. Grundlehren der Mathematischen Wissenschafte, 268. Springer, Heidelberg, 2011. 


\bibitem{bakerjensen} M.~Baker, D.~Jensen, Degeneration of linear series from the tropical point of view and applications, preprint, \texttt{arXiv:1504.05544}.

\bibitem{bpr} M.~Baker, S.~Payne, J.~Rabinoff, Nonarchimedean geometry, tropicalization, and metrics on curves, Algebraic Geometry 3 (2016), no. 1, 63--105. 

\bibitem{berk1} V.~Berkovich, Spectral theory and analytic geometry over non-Archimedean fields, Mathematical Surveys and Monographs, 33, American Mathematical Society, Providence, RI, 1990.


\bibitem{bosch-lutkebohmert} S.~Bosch, W.~L\"utkebohmert, 
Stable reduction and uniformization of abelian varieties. I., Math.~Ann.~270 (1985), no. 3, 349--379. 

\bibitem{bmv} S.~Brannetti, M.~Melo, F.~Viviani,  On the tropical Torelli map,
Advances in Mathematics 226 (2011) 2546--2586.


\bibitem{capsurvey} L.~Caporaso, Algebraic and tropical curves: comparing their moduli spaces,
Handbook of Moduli, Volume I, G. Farkas, I. Morrison (eds.), Advanced Lectures in Mathematics, Volume XXIV (2013), 119--160.

\bibitem{cha1} M.~Chan, Combinatorics of the tropical Torelli map, Algebra Number Theory 6 (2012), no.~6, 1133--1169.


\bibitem{cha14} M.~Chan, Topology of the tropical moduli spaces $M_{2,n}$, preprint, \texttt{arxiv:1507.03878}

\bibitem{cgp16} M.~Chan, S.~Galatius, S.~Payne, The tropicalization of the moduli space of curves II: topology and applications, preprint, \texttt{arxiv:1604.03176}

\bibitem{cv} M.~Culler, K.~Vogtmann, Moduli of graphs and automorphisms of free groups, Invent. Math. 84 (1986), no. 1, 91--119.

\bibitem{danilov} V.~I.~Danilov, Polyhedra of schemes and algebraic varieties, Mat.~Sb.~(N.S.) 139 (1975), no. 1, 146--158, 160.

\bibitem{deligne3} P.~Deligne, Th\'eorie de Hodge II, III,  Inst. Hautes \'Etudes Sci. Publ. Math. No. 40 (1971) 5--57 and 44 (1974), 5--77.

\bibitem{deligne-mumford} P.~Deligne, D.~Mumford, The irreducibility of the space of curves of given genus, Inst.~Hautes \'Etudes Sci.~Publ.~Math.~No.~36 (1969) 75--109.
 

\bibitem{gathmann-markwig} A.~Gathmann, H.~Markwig, The numbers of tropical plane curves through points in general position, J. Reine Angew. Math. 602 (2007), 155--177. 


\bibitem{getzler} E. Getzler, Resolving mixed Hodge modules on configuration spaces, Duke Math.~J.~96 (1999), no. 1, 175--203.


\bibitem{gubler} W.~Gubler, A guide to tropicalizations, Algebraic and combinatorial aspects of tropical geometry, 125--189, Contemp. Math., 589, Amer. Math. Soc., Providence, RI, 2013.

\bibitem{grw} W.~Gubler, J.~Rabinoff, and A.~Werner, Skeletons and tropicalizations, Adv.~Math.~294 (2016)
150--215.

\bibitem{hacking08} P.~Hacking, The homology of tropical varieties, Collect. Math. 59 (2008), no. 3, 263--273.


\bibitem{harris-morrison} J.~Harris and I.~Morrison, Moduli of curves, Graduate Texts in Mathematics, 187.  Springer--Verlag, New York (1998), xiv+366 pp.


\bibitem{gfan} A.~N.~Jensen, Gfan, a software system for Gr\"obner fans and tropical varieties, available at \texttt{http://home.imf.au.dk/jensen/software/gfan/gfan.html}. 

\bibitem{kajiwara} T.~Kajiwara, Tropical toric geometry,  Toric topology, 197--207, Contemp. Math., 460, Amer. Math. Soc., Providence, RI, 2008.



\bibitem{kkmsd} G. Kempf, F. Knudsen, D. Mumford, and B. Saint-Donat, Toroidal embeddings. I, Lecture Notes in Mathematics, Vol. 339, Springer-Verlag, Berlin, 1973.


\bibitem{knudsen1} F.~F.~Knudsen and D.~Mumford, The projectivity of the moduli space of stable curves. I. Preliminaries on ``det'' and ``Div'', Math. Scand. 39 (1976), no. 1, 19--55.

\bibitem{knudsen2} F.~F.~Knudsen, The projectivity of the moduli space of stable curves. II. The stacks $\mathcal{M}_{g,n}$.  Math. Scand. 52 (1983), no. 2, 161--199.  

\bibitem{knudsen3} F.~F.~Knudsen, The projectivity of the moduli space of stable curves. III. The line bundles on $\mathcal{M}_{g,n}$, and a proof of the projectivity of $\ov{M}_{g,n}$ in characteristic $0$.
Math. Scand. 52 (1983), no. 2, 200--212. 

\bibitem{ms} D.~Maclagan, B.~Sturmfels, Introduction to Tropical Geometry, Graduate Studies in Mathematics, 161, AMS, Providence, RI, 2014. vii+359pp.

\bibitem{markwig-thesis} H.~Markwig, The enumeration of plane tropical curves, Ph.D~dissertation, Technischen Universit\"at Kaiserslautern, 2006.

\bibitem{mik06} G.~Mikhalkin, Tropical geometry and its applications, International Congress of Mathematicians. Vol.~II, 827--852, Eur.~Math.~Soc., Z\"urich, 2006.

\bibitem{mz} G.~Mikhalkin, I.~Zharkov,  Tropical curves, their Jacobians and theta functions, Contemporary Mathematics 465, Proceedings of the International Conference on Curves and Abelian Varieties in honor of Roy Smith's 65th birthday (2007), 203--231.


\bibitem{pay} S.~Payne, Analytification is the limit of all tropicalizations, Math.~Res.~Lett.~16 (2009), no.~3, 543--556.



\bibitem{thuillier} A.~Thuillier, G\'eom\'etrie toro\"idale et g\'eom\'etrie analytique non archim\'edienne. Application au type d'homotopie de certains sch\'emas formels,
Manuscripta Math. 123 (2007), no. 4, 381--451.


\bibitem{tyomkin} I.~Tyomkin, Tropical geometry and correspondence theorems via toric stacks, Math. Ann. 353 (2012), no. 3, 945--995.

\bibitem{ulirsch-thesis} M.~Ulirsch, Tropical geometry of logarithmic schemes, Ph.D.~dissertation, Brown University, 2015.



\bibitem{viviani} F.~Viviani, Tropicalizing vs.~compactifying the Torelli morphism, Tropical and non-Archimedean geometry, 181--210, Contemp. Math., 605, Amer. Math. Soc., Providence, RI, 2013. 



\end{thebibliography}
\end{document}